# Contractible spaces and coalescent homotopies

Eduardo Francisco Rêgo

**Abstract**

This paper deals with the existence, or absence, of *coalescent contractions* of contractible spaces. These are the contractions such that when the tracks of any two points meet, at time $t_0$, they remain together thereafter. If a finite simplical complex $K$ is *collapsible*, then any *collapse* of $K$ encodes coalescent contractions of $K$. Examples of contractible spaces where no coalescent contractions exist are the *dunce hat* and *Bing's house*. We establish a criteria for contractible finite simplicial complexes that ensures there are no coalescent contractions: the *star-disc property*.

## 1  Introduction

This paper is about *contractible spaces* and the existence, or absence, of *coalescent* contractions. We say that a homotopy, $H : X \times I \longrightarrow Y$, is *coalescent* if

$$H(a, t_0) = H(b, t_0) \Rightarrow H(a, t) = H(b, t) \ \forall t \geq t_0$$

that is, when the tracks of any two points, $a, b \in X$ meet, at time $t_0$, they remain together thereafter. Another meaning, closely related to this one, is the notion of *opening time* of a contraction $H : Y \times I \longrightarrow Y$; this is the first time the slice functions $H_t$ stop beeing surjective: $t_{op} = \sup \{t : Y_s = Y \ , \ \forall s \leq t\}$, where $Y_s = H_s(Y)$.

In section 2 we revisit the *dunce hat*: this is the simplest example of a contractible finite simplicial complex which is non-collapsible, [1]. It is not collapsible, simply because there is no free simplex to start a collapse. The easiest proof that it is contractible is like this: one readily sees that the dunce hat can be embedded in euclidian space, $\mathbb{E}^3$, where it is a (strong) deformation retract of a neighbourhood (a regular neighbourhood) homeomorphic to the 3-ball, $B^3$, and is therefore contractible [1]. While it is reasonably easy to visualize, in $\mathbb{E}^3$, the embedding of the dunce hat, its thickening to a 3-ball and the corresponding deformation, it is far more difficult to visualize an actual contraction of the dunce hat, at least until one gets used to that mental exercise. What we shall do, in section 2, is to describe, in intrinsic terms, a contraction of the dunce hat; the main purpose is to understand there are contractions of the dunce hat with opening time $t_{op} = 0$ and have a clear picture of the neighbourhoods of the vertex.

In section 3 we establish the relation between coalescente contractions and opening time: if there is a coalescente contraction there is one with opening time 0.

In section 4 we prove that there are no coalescente contractions of the *Bing's house with two rooms*, another well known contractible and non-collapsible complex. In the processe we establish a general condition for finite contractible simplicial complexes that prevents the existence of any coalescent contractions: the *star-disk property*. The arguments used can be generally described as variations on Brower's fixed point theorem and degrees of maps between spheres;

In section 5 we prove that the dunce hat doesn't have any coalescent contractions. This comes in a first place as a consequence of the dunce hat not having *free* simplices, which also prevents it from being collapsible; of course, for a collapsible space, $K$, we can associate to each collapse, $K \searrow *$, a coalescent contraction. The dunce hat does not satisfy the star-disk property at its vertex, therefore some ad hoc arguments have to be trown in to complete the proof.

---

[1]For the background notions just invoked (simplicial complexes, collapsible complex, regular neigbhourhoods, etc.), see [2, Chapter 2] for a quick introduction and [3], [4] for deeper treatments.

## 2   From inside the dunce hat I: contraction

Recall that the dunce hat, $D$, is the quotient space obtained from a triangle (2-simplex) by identifying the three oriented edges where the orientations of two of the edges are pointing to the same vertex. See Figure 1, bellow.

The three vertices of the triangle are all identified to the same point $V$ of $D$, which we call the *vertex* of the dunce hat: as we shall see, this point plays an essential part in several properties of $D$. Represented in shade is a typical neighbourhood of $V$: it is easy to see that it is homeomorphic to a disc - the bottom right triangle $AVC$ - with two cones glued along generators, $VA$ and $VB$: see Figure 2 for a picture of an embedding of that neighbourhood in three space.

We start now, describing a contraction of the dunce hat. We look first at the cone corresponding to the top triangle in Figure 1. We divide it into three regions, depicted in Figures 1 and 2, as follows. Consider the perpendicular segment from $V$ to the foot $V_0$, and divide it in thirds by the points $V_1$, $V_2$. Consider next the two pairs of lines from $A$ to $V_1$ and $V_2$: each gives a loop based in $A$ and going through $V_i$. The perpendicular from a generic point $P$, in the generator $AV$, to its foot $P_0$, gets divided into thirds by the intersections with those loops, $P_1$ and $P_2$: one such perpendicular is represented in Figure 2.

We construct the contraction of $D$ as a sequence of several homotopy movements which we describe in order.

The first movement fixes all points outside the cone, and also those on the base of the cone, that is on loop $AV_0A$, as well as on the generator $AV$ and moves just those points in the interior of the top triangle, in Figure 1. The tracks of the points in the interior of the triangle follow the perpendicular segments: for each perpendicular segment $PP_0$, the homotopy fixes the end points, shrinks the segment $P_2P$ to the end point $P$ and stretches the segment $P_0P_2$ to a homeomorphism into $P_0P$. At the end of the movement, for any point $V'$ in the perpendicular segment $V_2V$ each of the two arcs from $A$ to $V'$ is sent homeomorphically to $AV$: so that the region between the loops $AV_0A$ and $AV_2A$ is stretched to give the original cone by identification of the two arcs $AV_2$ to $AV$, while the complementary region of the cone closes down to $AV$ - as a fan made of those loops $AV'A$.

The second part of the homotopy keeps moving just those points interior to the "fan region" just mentioned, and the movement takes place entirely within the segment $AV$, where they all lie by the end of the first movement. Recall that for a generic $P \in AV$, the "perpendicular segment" $P_2P$ was squeezed to $P$; now, the points $P_1$ move along $AV$, from $P$ to $A$ and correspondingly the two segments $PP_1$ and $P_2P_1$ expand again, keeping the other end points $P$ and $P_2$ fixed at $P$, to fill the segment $PA$ homeomorphically. Note that at the end of this process, the loop $AV_1A$ is shrunk to the point $A$. In Figure 3 bellow, we try to give a visual impression of the combination of the first two homotopy movements, just described. But note that it is not a realistic picture: the two arrows to the right, pointing downward, suggest rightly the first movement, that of closing the fan and the third arrow, pointing to the left, suggests the second movement, with the thick loop $AV_1A$ shrinking in the direction of point $A$; but in this respect it is a fake picture, since the second movement happens all inside the segment $VA$: lets say it's an "infinitesimal vision" from inside $AV$, allowing us to keep seeing all the elements of the closed fan - like the points $P_2$ and $P$ or the segments $P_2P_1$ and $P_2P_1$ in the picture - set apart ! Anyway, if the reader doesn't find this picture helpful he/she may simply ignore it.

The third part of the homotopy starts by pushing the segment $AV$, keeping the end points fixed, slightly inside the disc that corresponds to the bottom right triangle $AVB$, sending it to an inside arc $AdV$, as represented in Figure 2: by 'pushing' we mean to perform a strong deformation retraction of the disc through the shaded area. In the process, we must drag along all the points from the interior of the top most region of the triangle - the part of the fan above loop $AV_1A$ - that by the end of the first two movements were all laid inside the segment $AV$: the important point now is that since the loop $AV_1A$ was shrunk to the point $A$ that now remains fixed, the dragging gives a continuous function!

This is the stage where the image of the dunce-hat by the homotopy is no longer surjective. Figure 4 represents the state of affairs at this point.

The image of the dunce-hat $D$, after the small deformation retraction of disc $AVB$ from $AV$ across the disc $\Delta AVd$, is the image after the first two movements of the shaded region plus the segment $AV$ which remains the image, after those previous movements, of the upper region enclosed by the loop $AV_2A$ - recall that each one of the arcs $AV_2$ was mapped homeomorphically into $AV$ - and the loop $AV_1A$ which was shrunk to $A$; the image of the top region is now the arc $AdV$, with the segment $V_1V$ - represented thick - mapped homeomorphically from $A$ to $V$.

The interior of the triangle $\Delta AdV$ is what is missing from the image of $D$ at this point.

The fourth part of the contraction of $D$ is the continuation of the previous small deformation retraction all the way across the interior of the shaded region, from the arc $AdV$ to the complementary part of the boundary of that region, as suggested in Figure 5 bellow, composed with the final map of the first two movements. At the end of this deformation, the image of $D$ will just be the image of segment $VA$ plus the image of arc $AV_1A$ by that map. That map was the identity on $VA$ and on arc $AV_1A$ was the constant map to $A$. Therefore the image of $D$ at the end of this fourth part will be an arc, namely $VA$.

Since an arc is a contractible space, the fifth and last part of the contraction of $D$ will be a contraction of the arc $VA$: note that if we choose this final contraction to be a strong deformation retraction to the end point $V$, this point - the vertex of $D$ - will have remained fixed through out all the contraction of $D$.

This finishes the description of a contraction of the dunce hat, as seen from the inside.

We should note that many points that meet at a certain stage of the contraction - in the terminology of the introduction, they *coalesce* - separate again at later stages. As we shall see in the next section, that phenomenon can't be avoided in any contraction of $D$.

**Note**: In the first two parts of the contraction we only moved points in the cone with generator $AV$, until we were able to dig a hole near the vertex $V$, by pushing along $AV$; note that we can take the cone - equivalently the generator $AV$ - to be arbitrarily small. If we consider an *infinitesimal* cone (more formally, the uniform limit of a sequence of homotopies defined, as the one before, on a sequence of cones converging to vertex $V$...) the hole near $V$, that is non-surjectiveness, appears right from the start of the homotopy: $\forall \varepsilon > 0, \exists t < \varepsilon : H_t(D) \neq D$, where $H_t, t \in [0,1]$ denotes the contraction. That is not always the case for contractible spaces; the reader may convince himself of that fact after solving the following exercise on another famous contractible non-collapsible space, *Bing's house with two rooms* (Figure 6)

The house is built in the following way. Start with a closed rectangular box, and add a slab in the middle; then add two more boxes, as shown in Figure 6, one from the mid slab to the top and another one from the mid slab to the bottom; remove then the base and top from each of these two chimneys. You got a house with two separated rooms: you can get into the bottom room through the top entrance and into the top room from the bottom entrance (see[5]). Finally add two walls, one in each room as shown shaded in Figure 6. You can easily see that a regular neighbourhood of this object in $\mathbb{R}^3$ is a 3-ball, therefore it is contractible; but it is not collapsible since there are no free simplices.

**Exercise 1** *Describe a contraction of Bing's house: the essential point is the art of digging a hole somewhere, as we did for the dunce-hat...*

*Hint*: *having experienced with the dunce hat $D$ before, what you should expect now is to be able to dig a hole near one of the "special" corners of the house where, as in the vertex $V$ of $D$, there is some complexity in the way the several walls meet at that point (intuitively, and as we shall see later, one can not do such a thing as digging a hole near a manifold point); but, after some trials, you will soon find that you can not act locally near just one such point, instead you will have to work in a larger region with several of them.*

## 3 Coalescence and opening time

We have already said, in the introduction, what it means for a homotopy $H : X \times I \longrightarrow Y$ to be *coalescent*:

$$H(a, t_0) = H(b, t_0) \Rightarrow H(a, t) = H(b, t), \forall t \geq t_0$$

that is, if the the tracks of two points, $a$ and $b$, meet at time $t_0$, they remain together thereafter.

What we shall establish, in this section, is a specially nice feature of coalescent contractions.

From now on we will assume, unless otherwise stated, that all spaces are finite simplicial complexes, although most of the arguments can be carried over to more general classes of spaces, especially compact metric spaces.

Let $Y$ be a contractible space and $H : Y \times I \longrightarrow Y$ a contraction that is coalescent. Let, for each $t \in I$, $Y_t$ be the image of $Y$ at time $t$: $Y_t = H_t(Y) = H(Y \times t)$. Since $H$ is continuous, each $Y_t$ is a compact subspace of $Y$, homeomorphic to the identification space $Y/R$ where $R$ is the relation of coalescence at time $t$, $aRb \Leftrightarrow H_t(a) = H_t(b)$. The nice feature of coalescent contractions we were alluding to is that given any $Y_{t_0}$, what we see during the contraction, from time $t_0$ until the end, is also a contraction of subspace $Y_{t_0}$ in the ambient space $Y$. This is an immediate consequence of transgression for *identification spaces* that we recall:

**Theorem 2** *Let $f : X \longrightarrow Z$ be a continuous map and $p : X \longrightarrow Y$ an identification map such that $h = fp^{-1}$ is well defined (that is, $f$ is constant on each fibre $p^{-1}(y)$). Then:*

**a)** *$h : Y \longrightarrow Z$ is a continuous map.*

**b)** *In the particular case when $Y = X_f = X/K_f$ with $K_f$ the relation induced by $f$, defined by $aK_f b \Leftrightarrow f(a) = f(b)$, $h$ is injective. Therefore each continuous map, $f$, factors through the composition of a continuous surjection and a continuous injection.*

**c)** *If, furthermore, $f$ is surjective, then $h = fp^{-1} : X_f \longrightarrow Z$ is a continuous bijection and is a homeomorphism, $h : X_f \cong Z$, if and only if $f$ is an identification map.*

**Proof.** See, for instance, [6, Chapter VI]. ∎

We can now give the precise formulation of that nice feature we described above:

**Remark 3** *Let $Y_{t_0}$ be the image of $Y$ at time $t_0$, and let $p : Y \times [t_0, 1] \longrightarrow Y_{t_0} \times [t_0, 1]$ be defined by $p(y, t) = (H_{t_0} y, t)$. Since we are assuming that $H$ is* coalescent, *$J = p^{-1}H$ is single-valued, thus defines a function $J : Y_{t_0} \times [t_0, 1] \longrightarrow Y$. Because $Y \times [t_0, 1]$ is compact, $p$ is an identification map and so, by a) of the previous theorem, it is continuous. Furthermore, it's clear that for each $y \in Y_{t_0}$, $J(y, t_0) = y$, that is $J_{t_0}$ is the identity map on $Y_{t_0}$, and since $H_1$ is constant we have that $J_1$ is also constant and therefore $J$ realizes a contraction of $Y_{t_0}$ in the ambient space $Y$. Clearly, $J$ is also coalescent.*

Let $H : Y \times I \longrightarrow Y$ be a contraction. Define the *opening time*, $t_{op}$, by
$t_{op} = \sup \{t : Y_s = Y, \forall s \leq t\}$. The opening time is the time when space $Y$ first opens, in the sense that we start seeing non-surjectiveness. It is an easy exercise to show that $Y_{t_{op}} = Y$.

After our construction, in the previous section, of a contraction for the dunce-hat $D$, we explained in a note how it could be modified, through a limiting process, to yield a contraction $H$ such that $\forall \varepsilon > 0, \exists t < \varepsilon : H_t(D) \neq D$. This means that $t_{op} = 0$, the opening time is 0. The next result shows that if there is a coalescent contraction, there is a contraction with opening time 0.

**Theorem 4** *Let $H : Y \times I \longrightarrow Y$ be a coalescent contraction. Then there is a (coalescent) contraction $K_s$, $s \in I$, of $Y$ with opening time $s_{op} = 0$.*
**Proof.** *Let $t_{op}$ be the opening time of $H$ and consider $Y_{t_{op}}$. As we've just stated above, $Y_{t_{op}} = Y$: with the usual notation for $\varepsilon > 0$ neighbourhoods of subspaces $A \subset Y$, $N_\varepsilon(A) = \{y \in Y : d(y, A) < \varepsilon\}$, where $d$ is a metric for $Y$ and $d(y, A) = \inf \{d(y, a) : a \in A\}$, we have: given an arbitrary $Y_{t_0}$ and neighbourhood $N_\varepsilon(Y_{t_0})$, using the uniform continuity of $H$, due to the domain being compact, we can get a $\delta > 0$ such that $|t - t_0| < \delta \Rightarrow Y_t \subset N_\varepsilon(Y_{t_0})$; in particular, if $Y_{t_0} \neq Y$, since $Y_{t_0}$ is compact, we have $d(a, Y_{t_0}) > 0$ for $a \in Y - Y_{t_0}$ and so taking $\varepsilon < d(a, Y_{t_0})$ we have that $Y_t \neq Y$ for all $t$ such that $|t - t_0| < \delta$; therefore, by the definition of opening time, it has to be the case that $Y_{t_{op}} = Y$.*

*Refer back to the discussion in the remark above: with $p : Y \times [t_{op}, 1] \longrightarrow Y_{t_{op}} \times [t_0, 1] = Y \times [t_{op}, 1]$, we have that $J = Hp^{-1} : Y \times [t_{op}, 1] \longrightarrow Y$ realizes a (coalescent) contraction of $Y$ and clearly, by definition of $t_{op}$ as a supremum, $\forall \varepsilon > 0, \exists t < t_{op} + \varepsilon : J_t(Y) \neq Y$. The final contraction $K_s$, $s \in I$, is obtained by an obvious reparametrization of $J$: $t = s(1 - t_{op}) + t_{op}$; clearly, $s_{op} = 0$.* ∎

## 4 The star-disk property

We are now in the position of establishing the non-existence of coalescent contractions of the dunce-hat $D$. But we look first at Bing's house; let's denote it $B$. Let $H$ be a contraction of $B$ and assume that it is coalescent: by the previous theorem we can assume, without loss of generality, that the opening time for $H$ is 0, that is $\forall \varepsilon > 0, \exists t < \varepsilon : H_t(B) = B_t \neq B$. Let $t_n$, $n \in \mathbb{N}$, be a sequence convergent to 0, such that $B_{t_n} \neq B$ for all $n$, and, for each $n$, choose some point $c_n \in B - B_{t_n}$; without loss of generality, considering a subsequence if necessary, we can assume that $c_n$ converges to a point, say $c_n \longrightarrow c$. So, as soon as the contraction of $B$ gets started, we see points arbitrarily close to $c$ that are missed at arbitrarily early stages.

Assume $B$ is endowed with a triangulation by 2-simplices linearly embedded in 3-space: this means that the corners of the house are vertices of the triangulation and that the edges of the house are contained in the 1-skeleton. Looking at the *star* of point $c$, $St(c)$, the union of all simplices that contain $c$, we can assume without loss of generality that the points $c_n$ all belong to the same simplex $\sigma \in St(c)$. It should be clear, looking at $B$, that we can always get an embedded disc $\Delta$, made up of simplices of $St(c)$ including $\sigma$, and with $c$ in its interior. Here interior means interior of the manifold - the points not in boundary circle $\Sigma = \partial \Delta$ - not the interior of $\Delta$ as a subspace of $B$: for those points of $B$ which are not manifold points - there are eight of them at the corners of the two shaded walls in Figure 6 - we have several ways of constructing $\Delta$. We can also assume that, alongside with $c$, all the $c_n$ are also in the interior of $\Delta$. Now, we just need to look at (the first stages of) the contraction of disc $\Delta$ in $B$, that is, to look at the restriction of $H$ to $\Delta \times I$.

Let us first recall some results about self-maps of n-balls and spheres. Given the n-ball $B^n = \{x \in \mathbb{R}^n : \|x\| \leq 1\}$ and an interior point $x$, let $r_x : B^n - \{x\} \longrightarrow S^{n-1}$ denote the usual retraction given by *radial projection* from $x$.

**Lemma 5** *Let $x_n$, $n \geq 1$, be a sequence of interior points of $B^n$ converging to an interior point $x_0$ and $N$ a closed neighbourhood of $S^{n-1}$ disjoint from the set of points $x_i, i \geq 0$. Then the restrictions to $N$ of the set of retractions $\{r_{x_n}\}_{n \geq 0}$ is* equicontinuous, *that is,*

$$\forall \varepsilon > 0, \exists \delta > 0 : \forall x \in N, \forall n \geq 0, \|r_{x_n}(x) - x\| < \varepsilon$$

**Proof.** *We leave the proof as an exercise for the reader.* ∎

The next result says that continuous self-maps of the sphere which are sufficiently close to the identity are homotopic to the identity. In fact more can be said:

**Lemma 6** *Let $f : S^{n-1} \longrightarrow S^{n-1}$ be a continuous map such that no point is sent to its antipodal, that is, $f(x) \neq -x, \forall x \in S^{n-1}$. Then $f$ is homotopic to the identity.*
**Proof.** *Consider for each point $x \in S^{n-1}$ the segment in $B^n$ with end points $x$ and $f(x)$. We leave it for the reader to check that we obtain a homotopy $H : S^{n-1} \times I \longrightarrow S^{n-1}$ between the identity and $f$ by sliding from $x$ to $f(x)$ along those segments and compose with the radial projection, $r$, from the origin - more specifically $H(x,t) = r(x + t(f(x) - x))$ (the track of each point is the (shortest) arc of great circle that joins it to its image).* ∎

Consider again the disc $\Delta$ and the contraction $H$. Through some specific homeomorphism $(\Delta, \Sigma) \longrightarrow (B^2, S^1)$ we get the analogous of the two previous lemmas for the pair $(\Delta, \Sigma)$, namely a set of retractions $r_{c_n} : \Delta - \{c_n\} \longrightarrow \Sigma$, $n \geq 0$ (with $c_0 = c$), which is *equicontinuous* on any closed neighbourhood $N(\Sigma)$ disjoint from the $c_i, i \geq 0$, and a positive constant $k$ such that any continuous map of $f : \Sigma \longrightarrow \Sigma$ with $d(f, id_\Sigma) \leq k$
($d$ the sup metric) is homotopic to $id_\Sigma$. Fix a $\alpha$-closed neighbourhood of $\Sigma$ in $\Delta$, $N_\alpha(\Sigma) = \{x \in \Delta : d(x, \Sigma) \leq \alpha\}$ which doesn't contain any of the points $c_i$: by equicontinuity, there exists $\delta > 0$ such that $\forall x \in N_\alpha(\Sigma), \forall n \geq 0, \|r_{c_n}(x) - x\| < k/2$. Let, for each $\varepsilon > 0$, $M_\varepsilon$ denote the $\varepsilon$-neighbourhood of $\Delta$ in $B$: it is clear that for sufficiently small $\varepsilon$ there is a (strong deformation) retraction $r^\varepsilon$ of $M_\varepsilon$ into $\Delta$, such that, for all $y \in M_\varepsilon - \Delta$, $r^\varepsilon(y) \in \Sigma$ and $d(y, r^\varepsilon(y)) \leq \varepsilon$. Choose such an $\varepsilon$, with $\varepsilon \leq \min\{\alpha, \delta, k/2\}$. By uniform continuity of $H$, there is a time $t_0 > 0$ such that for all $x \in B$ and for all $t \leq t_0$ we have $d(x, H(x,t)) < \varepsilon$; choose a time $t_{c_n} < t_0$, call it $s$ to simplify notation and let $f$ be the restriction of $H_s$ to the disc $\Delta$. Let $g : \Delta \longrightarrow \Sigma$ be defined by $g = r_{c_n} \circ r^\varepsilon \circ f$: by the choice of $s = t_{c_n} < t_0$, $f(\Delta) \subset M_\varepsilon$, so $r^\varepsilon \circ f$ is well defined. Let $x \in \Sigma$: if $f(x) \in M_\varepsilon - \Delta$ we have $g(x) = r^\varepsilon \circ f(x)$ and $d(f(x), r^\varepsilon f(x)) \leq \varepsilon$, so $d(x, g(x)) \leq d(x, f(x)) + d(f(x), r^\varepsilon \circ f(x)) \leq \varepsilon + \varepsilon \leq k$; if $f(x) \in \Delta$, $g(x) = r_{c_n} \circ f(x)$, and since $d(x, f(x)) < \varepsilon \leq \min\{\alpha, \delta, k/2\}$, we have that $f(x) \in N_\alpha$; therefore, since $d(x, f(x)) < \delta \Rightarrow d(f(x), r_{c_n} \circ f(x)) < k/2$ we have $d(x, g(x)) \leq d(x, f(x)) + d(f(x), r_{c_n} \circ f(x)) \leq k/2 + k/2 = k$. In both cases we conclude that $\forall x \in \Sigma$, $d(g(x), x) < k$ and so, by the choice of $K$ above, the restriction of $g$ to $\Sigma$ is homotopic to the identity and so has *degree* 1; on the other hand, being the restriction of a map defined on the disc $\Delta$ it is nulhomotopic and so has *degree* 0, contradiction[2].

We have thus finished a proof that no contraction of space $B$ has opening time 0 - a result we hinted at in the note of the previous section - and as a corollary that no contraction of $B$ is coalescent:

**Theorem 7** *No contraction of Bing's house, $B$, has* opening time $t_{op} = 0$; *therefore no contraction of $B$ is coalescent.*

We should stress what was the key factor in the previous proof: given any point $c$ where space $B$ "opens up" (at time 0) and whatever simplex, $\sigma \in St(c)$, contains points arbitrarily close to $c$ from the complements of images $H_t(B)$, it is possible to get an embedded disc $\Delta$ containing $\sigma$ and with $c$ an interior point.

**Definition 8** *We say a point $c$ in a simplicial complex $C$ has the* star-disc *property if for each simplex $\sigma \in St(c)$ there is an embedded disc $\Delta$ that contains $\sigma$ and with $c$ an interior point (note that, in particular, this property implies that there are no free simplices and so $C$ is non-collapsible)*

It is straightforward - we leave it as an exercise for the reader - to generalize the proof to higher dimensions to get:

**Theorem 9** *Let $C$ be a finite, contractible simplicial complex such that every point $c$ has the* star-disc *property (therefore $C$ is non-collapsible)*
*Then, every contraction of $C$ has opening time $t_{op} > 0$ and so is non-coalescent.*

## 5  From inside the dunce hat II: non-coalescence

The dunce-hat $D$ doesn't satisfy the hypothesis of the theorem: there is exactly one point which doesn't have the required property, the vertex $V$. Referring back to Figures 1 and 2, it is clear that a simplex in $St(V)$ that comes from the right-bottom part, like triangle $AVB$, can not be extended to an embedded disc with $V$ an interior point. This fits in well with our note in the previous section, where it was observed that our construction of a contraction of $D$ would have opening time 0 if it was done with an *infinitesimal* cone; on the other hand any contraction of $D$ with $t_{op} = 0$ will have to open up at the vertex $V$ since it's easily seen that any point other than $V$ has the star-disc property.

Although the dunce-hat has contractions with opening time 0, it doesn't have any coalescent contractions, but to settle this some further arguments are needed to deal with the special situation of the vertex. We start by establishing a lemma to the effect of getting rid of a certain type of possible wilderness of the track of vertex $V$ under a contraction $H$ of $D$: for technical reasons, that will become apparent in the course of the final proof, we do not want the track $H(\{V\} \times I)$ to cover the complements of images $H_t(D)$ that arise after opening time.

**Lemma 10** *Let $H$ be a (coalescent) contraction of $D$ and $p$ any manifold point of $D$ (that is, $p$ is any point in the interior of the identification triangle in Figure 1). Then there is a (coalescent) contraction, $J$, such that the track of $V$ under $J$ misses $p$.*
**Proof.** *Let $\gamma(t) = H(V,t)$. Let $B$ be a closed neighbourhood of $p$, homeomorphic to a disc (see Figure 7) and disjoint from $V$.*

*Let $u_o$ be any time such that $\gamma(u_o) = p$,*
$t_0 = \inf \{t \in I : \gamma([t, u_0]) \subset B\}$ *and* $s_0 = \sup \{t \in I : \gamma([u_0, t]) \subset B\}$; *this means that in the interval $[t_0, s_0]$ the track of $V$ stays in disc $B$. Note that Figure 7 is an over simplified picture, with $\gamma([t_0, s_0])$ represented by a very simple polygonal arc: the point is that it could be much complicated, even with $\gamma([t_0, s_0])$ filling the whole of $B$ ! If there is another time $u_1$ outside this interval such that $\gamma(u_1) = p$, we have another interval, $[t_1, s_1]$, similarly defined and disjoint from the first: it's clear that we can repeat this process until we get a finite number of disjoint closed intervals $[t_i, s_i], i = 0, ..., m$ whose union contains the pre-image $\gamma^{-1}(p)$. For each $i = 0, ..., m$ we modify the homotopy $H$ in each time period $[t_i, s_i]$. Starting with $[t_0, s_0]$:*

---

[2]The reader who hasn't yet learned about this most important notion of *degree* - and the associated results on maps between spheres - introduced by Brower, is urged to do so: see [6, Chapters XV-XVII] or books on algebraic topology, for instance [2], [7], [8].

let $E$ be another disc surrounding $B$ as shown in Figure 7; consider an *isotopy of $D$*, $K : D \times I \longrightarrow D$ which is the identity outside $E$ and on $E$ is the identity on the boundary and moves interior points so as to push disc $B$ to a final position $k(B)$, $k = K_1$, disjoint from point $p$; as suggested in Figure 7, we can further assume that the isotopy is also fixed in an neighbourhood of an arc in $\partial B$ that contains the two points $\gamma(t_0) = H_{t_0}(V)$ and $\gamma(s_0) = H_{s_0}(V)$: just expand the shaded portion of the ring, keeping the complementary white part fixed, while pushing disc $B$ along the corresponding boundary arc. Let $N_\varepsilon(\gamma(t_0))$ and $N_\varepsilon(\gamma(s_0))$ be neighbourhoods of those two points, contained in $E$ and fixed through the isotopy. Let $\delta_0 > 0$ be sufficiently small so that $[t_0 - \delta_0, s_0 + \delta_0]$ is still disjoint from the other intervals $[t_i, s_i]$ and such that $\gamma([t_0 - \delta_0, t_0]) \subset N_\varepsilon(\gamma(t_0))$, $\gamma([s_o, s_0 + \delta_0]) \subset N_\varepsilon(\gamma(s_0))$. Now we modify $H$ in the time period $[t_0 - \delta_0, s_0 + \delta_0]$: consider the reparametrizations $s(t) = (1/\delta_0)[t - (t_0 - \delta_0)]$, $t_0 - \delta_0 \leq t \leq t_0$, and $r(t) = (1/\delta_0)[t - s_0]$, $s_0 \leq t \leq s_0 + \delta_0$, and define $J : D \times I \longrightarrow D$ by

$$\begin{aligned}
J(x,t) &= H(x,t) \text{ if } t \leq t_0 - \delta_0 \\
J(x,t) &= K(H(x,t), s(t)) \text{ if } t_0 - \delta_0 \leq t \leq t_0 \\
J(x,t) &= k(H(x,t)) \text{ if } t_0 \leq t \leq s_0 \\
J(x,t) &= K(H(x,t), 1 - r(t)) \text{ if } s_0 \leq t \leq s_0 + \delta_0 \\
J(x,t) &= H(x,t) \text{ if } s_0 + \delta_0 \leq t
\end{aligned}$$

In short, between times $t_0 - \delta_0$ and $t_0$ we combine, through a suitably reparametrization, $H$ with the isotopy $K$, next from $t_0$ until $s_0$ we compose $H$ with the final map $k = K_1$ of that isotopy, and finally between times $s_0$ and $s_0 + \delta_0$ we combine $H$ with the time reversed isotopy. It should be clear that under this new homotopy, $J$, the track of $V$ misses $p$ during period $[t_0, s_0]$ and was unchanged before or after. We next repeat the process, in turn, for the other time periods $[t_1, s_1], [t_2, s_2], ..., [t_m, s_m]$. The proof ends by noting that, since each changing of $H$ was achieved by combining it with an *isotopy*, if $H$ is coalescent so $J$ will be. ∎

With the previous lemma in hand, we can now adjust the arguments we used in proving theorem 7 to the exceptional situation of vertex $V$ in $D$.

Let $H$ be a *coalescent* contraction of $D$. Choose any *manifold* point $p$ in a complement of an image $D - H_t(D)$, for some $t > t_{op}$ (of course we could assume, without loss of generality, that $t_{op} = 0$ but there is no gain in assuming that). By the previous lemma, there is a coalescent contraction $J$ such that the track of $V$ misses $p$. By usual reasoning with uniform continuity and compactness we can assume there is a neighbourhood of $V$, $N_\varepsilon(V)$, such that the track of each of its points also misses $p$, that is $p \notin J(N_\varepsilon(V) \times I)$. Refer back to Figures 1 and 2 that represent a neighbourhood of $V$: consider in each of the two cones a closed segment of a generator with end point $V$ and contained in $N_\varepsilon(V)$ and glue to the dunce-hat, by an homeomorphism, a disc $E$ along the union of both segments as represented in Figure 8.

Labelled $L$, the union of those two segments is represented thick. Denote the resulting space by $F$.

**Remark 11** *Observe that the only points $c \in F$ that do not satisfy the* star-disc *property (Definition 8) are the points on the* free arc *of $E$, including the two end points where that arc meets $L$; but even for these extreme points the property only fails for those simplices $\sigma \in St(c)$ in the complement of $D$.*

$F$ is of course contractible: there is an obvious strong deformation retract of $F$ into $D$ - across the disc $E$ from the free arc to the gluing arc $L$ - and therefore we can follow this first homotopy by a contraction of $D$. But we are going to construct a different contraction of $F$. Consider a regular neighbourhood $N$ of $L$ in $E$, homeomorphic to $L \times I$ (represented shaded gray in Figure 8). We define a contraction $K : F \times I \longrightarrow F$ in two steps, first for $t \leq 1/2$ and second for $1/2 \leq t \leq 1$. In the first part: if $x \in D$ we define $K(x,t) = J(x, 2t)$, if $x \in \overline{E - N}$ we have $K(x,t) = x$ and for those points in the region $N \cong L \times I$ the homotopy stretches each stalk $S_x = \{x\} \times I$, $x \in L$ (one of those is depicted in Figure 8), keeping the top end-point in $\overline{E - N}$ fixed, so as to follow $x$ along its track under $J$ - we leave it as an exercise for the reader to define suitable formulas for this. At the end of this first part, we have $K_{1/2}(F) = E \cup J(L \times I)$: $D$ was contracted to a point, say $q$, by $J$ (with double speed), and for each $x \in L$ we developed its former track under $J$ as we've just explained. In the second part of the homotopy, for times $1/2 \leq t \leq 1$, we contract $K_{1/2}(F)$ in itself, to the point $q$, in the obvious way: we first use a strong deformation retraction of $E \cup J(L \times I)$ into $J(L \times I)$, across the disc $E$, as alluded above; then we contract $J(L \times I)$ to $q$ erasing all the tracks - we also leave the appropriate formulation of this *erasing* to the reader.

It is easy to check that each step of the construction of $K$ preserves the given coalescence of $J$, so $K$ is also coalescent.

By our assumptions about $J$, and the fact that $L \subset N_\varepsilon(V)$, the tracks under $K$ of all points $x \in L$ miss $p$ and so $K_{1/2}(F) \neq F$; therefore the opening time for $K$ is strictly less than $1/2$ - recall that at its opening time any homotopy is still surjective. So, by theorem 4, we can assume that $t_{op} = 0$ - starting afresh at opening time and reparametrizing, which we do here by sending $[t_{op}, 1/2]$ to $[0, 1/2]$ and keeping the other times fixed. Now, as in the proof of theorem 7, let $t_n$, $n \in \mathbb{N}$, be a sequence of times convergent to 0, and $c_n \longrightarrow c$ a sequence of points, all in a simplex $\sigma \in St(c)$, and such that, for each $n$, $c_n \in F - K_{t_n}(F)$. Since for all $t \leq 1/2$ we have $E \subset K_t(F)$, the simplex $\sigma$ is not contained in $E$ and therefore - see Remark 11, above - the point $c$ has the star-disc property that allows us to carry on the argument of theorem 7 to arrive at a contradiction.

We have thus proved:

**Theorem 12** *Any contraction of the dunce-hat is non-coalescent.*

We end this final section by leaving two problems:

**Problem 13** *Generalize last theorem to arbitrary finite contractible simplicial complexes.*

The troublesome part of the generalization will be to work out appropriate analogues of the previous gluing of disc $E$: that is, to add appropriate pageants near a point where the star-disc property fails.

The second problem is just a fun-exercise about contractions: it asks about the existence of *monotone* contractions (definition in the statement) of a contractible space. The answer is yes in general[3], but I don't know the answer with the coalescence condition.

**Problem 14** *Let $H$ be a (coalescent) contraction of space $Y$. Is there always a (coalescent) contraction of $Y$ that is descending, or monotone, in the sense that $t' > t \Rightarrow H_{t'}(Y) \subset H_t(Y)$?*

# References

[1] E. C. Zeeman, On The Dunce Hat, Topology, vol.2 (1964), 341-358.

[2] C. R. F. Maunder, *Algebraic Topology*, Van Nostrand Reinhold Company, London (1972).

[3] C. Rourke and B. Sanderson, *Introduction to piecewise-linear topology*, Springer-Verlag (1972).

[4] T. B. Rushing, *Topological embeddings*, Academic Press (1973).

[5] R. H. Bing, *The Geometric Topology of 3- Manifolds*, Amer. Math. Soc. Colloquium Publications, vol. 40, Amer. Math. Soc. (1983).

[6] J. Dugundji, *Topology*, Allyn and Bacon, Inc., Boston (1966).

[7] J. F. Davis and P. Kirk, *Lecture Notes in Algebraic Topology*, Graduate Studies in Mathematics, vol. 35, Amer.Math.Soc. (2001).

[8] G. E. Bredon, *Topology and Geometry*, Springer-Verlag (1993).

Eduardo Francisco Rêgo
*Adress*: Largo 25 de Abril, 70 ; 4900-027 Afife Portugal
*E-mail address*: efrancisco.rego@gmail.com ; eerego@fc.up.pt
*Affiliation*: retired member of *Departamento de Matemática - FCUP Universidade do Porto*, Portugal.

---

[3]The idea of the proof is easy. Here is a very sketchy description. The homotopy is built in four steps. 1- Consider a top dimensional simplex $\sigma \in K$. Keeping its boundary and exterior fixed, expand its interior, in a increasing monotone fashion to cover the whole complex $K$. The inverse homotopy is then monotone decreasing. 2- Consider a strong deformation retraction of $K$ into $\sigma$. Such a deformation exists by standard results: see Dugundji. Note that up to this point the images of the homotopy slices are the whole complex. 3 - Reverse the homotopy built in the first step. The final image is the cimplex $\sigma$. 4- Finally contract $\sigma$ to one of its points (which can be done in a decreasing monote way)